\documentclass{amsart}
\usepackage[dvips]{epsfig}
\usepackage{graphicx}
\usepackage{latexsym}
\usepackage{amsmath}
\usepackage{amsthm}
\usepackage{amssymb}
\RequirePackage[colorlinks=true,citecolor=blue,urlcolor=blue]{hyperref}

\usepackage{subfigure}

\usepackage{amsfonts,multirow,rotating}
\usepackage{color}
\definecolor{gris}{gray}{0.4}
\usepackage{dsfont}

\usepackage{mathtools}

 \usepackage[applemac]{inputenc}

\newtheorem{thm}{Theorem}
\newtheorem{assumption}[thm]{Assumption}
\newtheorem{lem}[thm]{Lemma}

\newtheorem{cor}[thm]{Corollary}
\newtheorem{defi}[thm]{Definition}

\newtheorem{prop}[thm]{Proposition}
\newtheorem{rk}[thm]{Remark}

\newtheorem{algo}[thm]{Algorithm}

\newcommand{\vip}{\vskip.2cm}
\newcommand{\COMMENTAIRE}[1]{}

\newcommand{\field}[1]{\mathbb{#1}}
\newcommand{\EE}{\field{E}}

\newcommand{\GG}{\field{G}}

\newcommand{\NN}{\field{N}}
\newcommand{\PP}{\field{P}}

\newcommand{\RR}{\field{R}}
\newcommand{\TT}{\field{T}}

\newcommand{\Bb}{{\mathcal B}}

\newcommand{\Ff}{{\mathcal F}}

\newcommand{\Hh}{{\mathcal H}}
\newcommand{\Ll}{{\mathcal L}}

\newcommand{\Pp}{{\mathcal P}}
\newcommand{\Qq}{{\mathcal Q}}

\newcommand{\Ttransition}{\boldsymbol{\Pp}}

\begin{document}

\title[Bandwidth selection for kernel density estimation in a BMCM]{Local bandwidth selection for kernel density estimation in bifurcating Markov chain model}

\author{S. Val\`ere Bitseki Penda and Angelina Roche}

\address{S. Val\`ere Bitseki Penda, IMB, CNRS-UMR 5584, Universit\'e Bourgogne Franche-Comt\'e, 9 avenue Alain Savary, 21078 Dijon Cedex, France.}
\email{simeon-valere.bitseki-penda@u-bourgogne.fr}

\address{Angelina Roche, CEREMADE, CNRS-UMR 7534,
Universit\'e Paris-Dauphine, Place du mar\'echal de Lattre de Tassigny,
75775 Paris Cedex 16, France.}
\email{roche@ceremade.dauphine.Fr}

\begin{abstract}
We propose an adaptive estimator for the stationary distribution of a bifurcating Markov Chain on $\mathbb R^d$. Bifurcating Markov chains (BMC for short) are a class of stochastic processes indexed by regular binary trees. A kernel estimator is proposed whose bandwidth is selected by a method inspired by the works of Goldenshluger and Lepski \cite{GL11}. Drawing inspiration from dimension jump methods for model selection, we also provide an algorithm to select the best constant in the penalty. \end{abstract}

\maketitle

\textbf{Keywords}: Bandwidth selection ; Bifurcating autoregressive process ; Bifurcating Markov chains ; Binary trees ; Nonparametric kernel estimation.

\textbf{Mathematics Subject Classification (2010)}: 62G05, 62G10, 62G20, 60J80, 60F05, 60F10, 60J20, 92D25.


\section{\textsc{Introduction}}

First introduced by Basawa and Zhou (2004) \cite{BZ04}, Bifurcating Markov chains models (BMCM for short) has recently received particular attention for its application to cell lineage study. Guyon (2007) \cite{Guyon}, have proposed such a model to detect cellular aging in \emph{Escherichia Coli} and proved laws of large numbers and central limit theorem for this class of stochastic process. Bitseki-Penda~\emph{et~al.}~(2014)~\cite{BDG14} have then completed these asymptotic results and proved concentration inequalities.  

To the best of our knowledge, kernel density estimation for the BMCM were considered first  by Doumic \& \textit{al.}  \cite{DHKR1},  where they estimate the division rate of population of cells reproducing by symmetric division, \textit{i.e.} cell reproduction where each fission produces two equal daughter cells. After this work, Bitseki \& \textit{al.}  \cite{BHO}, have used the wavelets methodology to study the nonparametric estimation of the density of BMCM. They propose an adaptive estimator in dimension 1. Recently, Bitseki and Olivier \cite{BO16} have studied the Nadaraya-Watson type estimators of a BMCM that they called nonlinear bifurcating autoregressive process. The latter model can be seen as an adaptation of nonlinear autoregressive process on binary regular tree. We mention that, except in \cite{BHO}, all the estimations done in the previous works are non adaptive. In particular, the question of data-driven bandwidth selection was not addressed in \cite{DHKR1} and \cite{BO16}. The main objective of this work is then to propose a data-driven method  for choosing the bandwidth for the kernel estimator of the invariant measure in the multi-dimensional BMCM, following ideas from the works of Goldenshluger and Lepski (2011) \cite{GL11}.

The idea of the method is to select the bandwidth minimizing an empirical criterion imitating the bias-variance decomposition of the risk of the kernel estimator.  More precisely, let $\Hh$ be a  collection of bandwidths and let $(\widehat{\nu}_{h})_{h\in \Hh}$ be a family of kernel estimators of an unknown density $\nu$.  Then, we select the bandwidth $\widehat{h}$ as 
\begin{equation*}
\widehat{h} = \arg\min_{h\in\Hh}\{\widehat{A}(h) + bV(h)\}
\end{equation*}
with $A(h) $ an empirical version of the bias of the estimator $\widehat\nu_h$ and $V(h)$ a penalty term with the same order than the variance. 
The now so-called \emph{Goldenshuger-Lepski} methodology (\cite{GL11}, but also \cite{GL14,Rebelles15a, Rebelles15b,PR16,LM16}), initially developed for density estimation, has been applied in many contexts such as deconvolution problems \cite{CL13}, conditional cumulative distribution function estimation \cite{CR14,Chagny}, regression problems \cite{Chagny,CR16}, conditional density estimation \cite{Chagny, BLR16}, hazard rate estimation \cite{Chagny}, white noise model \cite{Lepski}, kernel empirical risk minimization (including robust regression) \cite{CL15}, L\'evy processes \cite{BL15}, Cox model \cite{GLT16}, stochastic differential model \cite{DG16}. Under suitable assumptions on the kernel, it is shown in \cite{GL11} that this selection rule leads to a minimax adaptive estimator on a general class of regular functions, for a general class of $\mathbb L_s$-risks. Pointwise versions of the Goldenshluger-Lepski selection rule have been less considered (Comte et Lacour (2013) \cite{CL13}, Rebelles (2015) \cite{Rebelles15b} and Chagny and Roche (2016) \cite{CR16}). The interest of such an approach is that the bandwidth is selected with a local criterion which realizes the best bias-variance compromise at the point where the estimator is calculated. On the contrary, integrated versions of the Goldenshluger-Lepski selection rule select the same bandwidth at all points. Let us mention that, in all the articles cited above, the theoretical results rely on concentration inequalities for sums of i.i.d. random variables such as Bernstein Inequality or Talagrand Inequality. In our context, such results are not applicable. Hence, we prove a Bernstein-type Inequality for functionals of BMC, where the functions are kernels and convolution of kernels. Compare to those obtained in~\cite{BHO}, our inequalities are more complete in the sense that the deviation parameter can take all the positive values. More precisely, its values do not depend on the size of the samples, with is essential for our theoretical results.

Ideally,  the penalty term $V$, called ``the majorant'' by Goldenshluger and Lepski, depends entirely on the kernel and the observations (it does not depend on the density $\nu$). The selected estimator is then $\widehat{\nu}_{\widehat{h}}$. However, for the BMCM, the variance term $V$ in the previous selection rule contains a term which may depend on the unknown density $\nu$. Moreover, this term is generally not estimable from observations. To resolve this problem, we propose a modification of Goldenshluger-Lepski rule's selection, inspired by the works of Lacour and Massart (2016) \cite{LM16}. As suggested in \cite{LM16}, the constant term in $bV(x,h)$ is then selected automatically from the data with an algorithm inspired by the works of Arlot and Massart (2009) \cite{AM09}.

The paper is organized as follows. The model is defined in Section~\ref{sec:definitions}. Section~\ref{sec:estimation} is devoted to the definition of the estimator. In Section~\ref{sec:simus}, we provide a numerical study of our estimator. The proofs are given in Section~\ref{sec:proofs}. 

\section{Definitions}
\label{sec:definitions}

We are now going to give a precise definition of a BMC. First we note that this class of stochastic processes has been introduced by Guyon \cite{Guyon} in order to understand the mechanisms of cell division. Indeed, these stochastic processes are well adapted to study a population (or more generally, any dynamic system) where each individual (or more generally, each particle) in one generation gives birth to two individuals in the next one. In the sequel, we will then use the language of the population dynamic to define the sets of interest.  

Let $(\Omega,\mathcal F, (\mathcal F_m, m\in\NN),\PP)$ be a filtered probability space. Let $(X_u, u\in \TT)$ be a sequence of random variables defined on $(\Omega, \PP)$, taking values in $\RR^{d}$, where $d\geq 1$, and indexed by the infinite binary tree
$\TT = \bigcup_{m=0}^\infty \{0,1\}^m$, with the convention that $\{0,1\}^{0} = \emptyset$. We equip $\RR^{d}$ with its usual Borel $\sigma$-field. Now, we will see $\TT$ as a given population. Then each individual $u$ of this population is represented by a sequence of $0$'s and $1$'s, and has two descendants, $u0$ and $u1$. The initial individual of the population is $\emptyset$. For all $m\in\NN$, let $\GG_{m}$ be the set of individuals belonging to the $m$-th generation, and $\TT_{m}$ the set of individuals belonging to the $m$ first generations. We have:

\begin{equation*}
\GG_{m}=\{0,1\}^m, \, \, \TT_{m} = \bigcup_{q=0}^{m} \GG_{q} \, \, \text{ and } \, \, \TT = \bigcup_{m\geq 0}\mathbb G_m. 
\end{equation*}
 
For an individual $u\in\GG_{m}$, we set $|u|:=m$ its length (i.e. the generation to which it belongs). 

\subsection{Bifurcating Markov chain \cite{Guyon}}

\begin{defi}[$\mathbb T$-transition probability] 

Let $\Pp: \RR^{d} \times \Bb((\RR^{d})^2) \to [0,1]$ (with $\Bb(\RR^{d})^2 = \Bb(\RR^{d}) \otimes \Bb(\RR^{d})$ the usual product $\sigma$-field on $(\RR^{d})^2$). Then $\Pp$ is a \emph{$\mathbb T$-transition probability} if 
\begin{itemize}
\item $x\to \Pp(x,A)$ is measurable for all $A\in\Bb(\RR^{d})^2$. 
\item $A\to \Pp(x,A)$ is a probability measure on $((\RR^{d})^2,\Bb(\RR^{d})^2)$ for all $x\in \RR^{d}$. 
\end{itemize} 
\end{defi}

For a $\Bb(\RR^{d})^{3}$-measurable function $f: (\RR^{d})^3 \to \RR$, we denote (when it is defined) by $\Pp f$ the $\Bb(\RR^{d})$-measurable function
 
\begin{equation*}
x\in \RR^{d} \mapsto \int_{(\RR^{d})^2} f(x,y,z) \Pp(x,dy,dz).
\end{equation*}

\begin{defi}[Bifurcating Markov chain]
Let $\mu$ be a probability measure on $(\RR^{d},\Bb(\RR^{d}))$ and $\Pp$ a $\mathbb T$-transition probability. We say that $(X_u, u\in \mathbb T)$ is a ($\mathcal F_m$)-bifurcating Markov chain with initial distribution $\mu$ and $\mathbb T$-transition probability $\Pp$ (denoted $\mathbb T$-BMC in the sequel) if 
\begin{itemize}
\item $X_u$ is $\mathcal F_{m}$ measurable for all $u\in\mathbb G_m$. 
\item $X_{\emptyset}$ has distribution $\mu$. 
\item  For all $m\in\NN$, and for all family $(f_u, u\in\mathbb G_m)$ of $\Bb(\RR^{d})$-measurable functions from $(\RR^{d})^3$ to $\RR$, 
$$\EE\left[\prod_{u\in\mathbb G_m}f_u(X_u,X_{u0},X_{u1})|\mathcal F_m\right]=\prod_{u\in\mathbb G_m}\Pp f_u(X_u),$$
where $u0:=(u,0)\in\mathbb G_{m+1}$ and $u1:=(u,1)\in\mathbb G_{m+1}$. 
\end{itemize}
\end{defi}

\subsection{Tagged-branched chain \cite{Guyon, BO16}}
Let $ (X_u, u\in \mathbb T)$  be a $\mathbb T$-BMC with initial distribution $\mu$ and $\mathbb T$-transition probability $\Pp$. We denote by $\Pp_0$ and $\Pp_1$ respectively the first and second marginals of $\Pp$. More precisely
 
\begin{equation*}
\Pp_0(x,B) = \Pp(x,B\times \RR^{d}) \, \, \text{ and } \, \, \Pp_1(x,B) = \Pp(x,\RR^{d}\times B),
\end{equation*}

for all $x\in \RR^{d}$ and all $B\in\Bb(\RR^{d})$. Let $\Qq$ be the mixture of $\Pp_0$ and $\Pp_1$ with equal weights

\begin{equation*}
\Qq = \frac12\Pp_0 + \frac12\Pp_1.
\end{equation*}

The Markov chain $Y:=(Y_m)_{m\in \NN}$ on $\RR^{d}$ with initial value $Y_0:=X_{\emptyset}$ and transition probability $\Qq$ is called the \emph{tagged-branch chain}.

In all the paper, we will denote by $\Qq^m$ the
$m$th iterated of $\Qq$ recursively defined by the formulas
\begin{equation*}\label{eq:transition}
\Qq^0(x,\cdot)=\delta_x \quad \text{and} \quad
\Qq^{m+1}(x,B)=\int_{S}\Qq(x,dy)\Qq^m(y,B) \quad \forall B\in
\mathcal{B}(\RR^{d}).
\end{equation*}
It is well known that $\Qq^m$ is a transition probability in
$(\RR^{d},\mathcal{B}(\RR^{d}))$. In particular, we have $\EE[f(Y_m)] = \mu\Qq^{m}f$ for all measurable function $f: \RR^{d} \to \RR$.

In the sequel, we will assume that the Markov chain $Y$ is ergodic- that is say, that there exists a unique distribution $\nu$ on $(\RR^{d},\Bb(\RR^{d})$ such that, for all measurable function $f: \RR^{d} \to \RR$,
 
$$
\lim_{m\to\infty}\EE[f(Y_m)]=\int_{\RR^{d}} f d\nu.
$$
We will also assume that the distribution $\nu$ has a density, that we also denote by $\nu$, with respect to the Lebesgue measure.
 
As previous works have shown (see for example \cite{Guyon, BDG14}), the analysis of a BMC $(X_{u}, u\in \TT)$ is strongly related to the asymptotic behavior of the tagged-branched chain $(Y_{m}, m\in \NN)$, and therefore to the knowledge of the invariant distribution $\nu$. We stress that this distribution is unknown and it is not directly observable, in such a way that its estimation from the data is of great interest.  The aim is to estimate $\nu$ from the observation of a subpopulation $(X_{u}, u\in\TT_{n})$.

\section{Estimation of the stationary distribution $\nu$}
\label{sec:estimation}
\subsection{Definition of the estimator}

We suppose that we observe the process $(X_u, u\in\mathbb T)$ up to the $n$-th generation. We denote by $|\mathbb T_n|=2^{n+1}-1$ the cardinality of $\mathbb T_n$. Based on the observation of $(X_u, u\in\mathbb T_n)$, we propose the following estimator of $\nu$

\begin{equation}\label{eq:nuhat}
\widehat \nu_h(x)=\frac1{|\mathbb T_n|}\sum_{u\in\mathbb T_n}K_h(x-X_u) \quad \forall x\in \RR^{d},
\end{equation}
where for all $x=(x_1,...,x_d)^t\in \RR^d$ and for all $h=(h_1,...,h_d)^t\in]0,+\infty[^d$
$$
K_h(x)=K(x_1/h_1,...,x_d/h_d)/\prod_{j=1}^dh_j
$$

and $K$ is a kernel, that is to say a function $K:\RR^d\to \RR$ which verifies $\int_{\RR^d}K(u)du=1$. 

The vector $h$ is usually called a bandwidth. In the sequel, $|h|$ will denote the product $h_{1}\times\cdots\times h_{d}$. It is well known in the kernel estimation theory that the choice of $h$ is of great interest. Indeed the amount of smoothing is controlled by a judicious choice of the bandwidth. Now, in order to tackle the issue of the choice of this bandwidth, we will need the following assumptions. 

\begin{assumption}[Assumptions on the kernel]\label{hyp:noyau}
$$
\|K\|_1=\int_{\RR^{d}} |K(t)|dt <+ \infty, \quad \|K\|_{2}^{2} = \int_{\RR^{d}} |K(t)|^{2}dt < +\infty \quad \text{and} \quad \|K\|_{\infty} = \sup_{t\in \RR^{d}} |K(t)|< +\infty. 
$$
\end{assumption}

\begin{assumption}[Uniform geometric ergodicity condition] \label{hyp:uniform_geometric_ergodicity}
There exists two constants $\rho \in (0,1/2)$ and $M>0$ such that for all bounded $\nu$-integrable function $g$, for all $x\in \RR^{d}$ and $m\geq 0$
\begin{equation*}\label{eq:lyapunov}
|\Qq^{m}g(x) - \int_{\RR^{d}} g d\nu|\leq M\|g\|_{\infty}\rho^{m}.
\end{equation*}
\end{assumption}

This assumption is verified for severals models of BMC. We refer for example to \cite{BO16} where Bitseki and Olivier have shown it for NBAR processes (see \cite{BO16} Lemma 20). For a precise definition of NBAR processes, we refer to section \ref{sec:simus}.

\begin{assumption}[Assumption on $\Pp$, $\Pp_{0}$, $\Pp_{1}$, $\Qq$ and $\nu$]\label{hyp:transition}
We assume that the transitions $\Pp$, $\Pp_{0}$, $\Pp_{1}$ and $\Qq$ admit densities with respect to the Lebesgue measure that we denote with the same notations. Moreover, we assume that
$$
\|\Pp\|_{\infty} < +\infty, \, \, \, \|\Pp_{0}\|_{\infty} < +\infty, \, \, \,  \|\Pp_{1}\|_{\infty} < +\infty, \, \, \, \|\nu\|_{\infty} < +\infty, \, \, \, \text{and} \, \, \, \|\Qq\|_{\infty} < +\infty.
$$ 
\end{assumption}
 
\subsection{Bias-variance decomposition}
We consider a pointwise quadratic risk 
$$\EE[(\widehat \nu_h(x)-\nu(x))^2],$$
where $x\in\RR^d$ is a given point. 

The results of \cite{BO16} (see the proof of Proposition 21) allows to obtain the following upper-bound on the risk, where we have make explicit the constants that appear. 

\begin{prop}
Under Assumption~\ref{hyp:noyau} to Assumption~\ref{hyp:transition}, we have 
\begin{equation}\label{eq:BV}
\EE[(\widehat \nu_h(x)-\nu(x))^2] \leq 2(K_h\ast \nu(x)-\nu(x))^2 + 2\frac{C(P,\nu)}{|\mathbb T_n||h|}, 
\end{equation}
where $\ast$ denotes the convolution product $f\ast g(x)=\int_{\RR^{d}} f(x-t)g(t)dt$ for all functions $f, g$ integrable over $\RR^{d}$ and $$C(P,\nu) = \frac{C_{I}}{(\sqrt{2}-1)^{2}}$$ with $$C_{I} = (1 + \frac{1}{1 - 2\rho^{2}})(\|\Qq\|_{\infty} + \|\nu\|_{\infty})^{2} + M^{2} + C_{P}$$ and $$C_{P} = 2 \|K\|_{2}^{2}(\|\Qq\|_{\infty} + \|\nu\|_{\infty}) + \|\Pp\|_{\infty} + \|\nu\|_{\infty}(\|K\|_{2}^{2} + \|\Pp_{0}\|_{\infty} + \|\Pp_{1}\|_{\infty})$$
\end{prop}

Inequality~\eqref{eq:BV} can be seen as a bias-variance decomposition in the sense that 
$$\EE_\nu[\widehat\nu_h(x)]=\frac1{|\mathbb T_n|}\sum_{u\in\mathbb T_n}\EE_\nu[K_h(x-X_u)]=\int_{\RR^{d}} K_h(x-t)d\nu(t)=K_h\ast\nu(x),$$
where $\EE_\nu$ is the expectation with respect to the measure $\nu$ and $C(P,\mu)/ |\mathbb T_n||h|$ is an upper-bound on the variance term $\EE\left[(\widehat\nu_h(x)-K_h\ast\nu(x))^2\right]$.

\subsection{Selection rule}
Given a family of bandwidths $\mathcal H_n\subset[0,+\infty[^d$, we have a family of estimators $(\widehat\nu_h)_{h\in\mathcal H_n}$ and the aim is to select an estimator in this family with risk close to the unknown oracle risk 
$$\mathbb E\left[(\widehat\nu_{h^*}(x)-\nu(x))^2\right] \text{ with } h^*={\arg\min}_{h\in\mathcal H_n}\mathbb E\left[(\widehat\nu_{h}(x)-\nu(x))^2\right].$$
Imitating the decomposition given in Equation~\eqref{eq:BV} we could consider the following selection rule inspired by the work of \cite{GL11}
\begin{equation}\label{eq:critere}
\widehat h = {\arg\min}_{h\in\mathcal H_n}\{A(x,h)+bV(x,h)\}
\end{equation}
where 
\begin{itemize}
\item $\mathcal H_n\subset[0,+\infty[^d$ is a finite collection of bandwidths;
\item $A(x,h)=\max_{h'\in\mathcal H_n}\left((\widehat\nu_{h'}(x)-K_h\ast\widehat\nu_{h'}(x))^2-aV(x,h')\right)_+$ with $b\geq a\geq 1$;
\item $V(x,h)= C(P,\mu)\log(|\mathbb T_n|)/  |\mathbb T_n||h|$.
\end{itemize}

The constants $a$ and $b$ can be different, as suggested in Section 5 of \cite{LM16}. In Section~\ref{sec:slope_estim}, we provide a method to choose both $a$ and $b$, as well as the quantity $C(P,\mu)$ appearing in the variance term $V(x,h)$. The lower bound $a=1$ is a critical minimal value, in the sense that the procedure fails for lower values of $a$ (see \cite{LM16}). 

We prove the following oracle-inequality on the selected estimator $\widehat\nu_{\widehat h}$. 
\begin{thm} \label{thm:oracle}Under Assumption~\ref{hyp:uniform_geometric_ergodicity}, if $\min_{h\in\mathcal H_n}|h|\geq \log(|\mathbb T_n|)/|\mathbb T_n|$, then there exists a minimal value $a_{\min}>0$ independent of $n$, such that, for all $a>a_{\min}$, 
\begin{equation}\label{eq:ineg_oracle_ponctuelle}
\EE\left[(\widehat\nu_{\widehat h}(x)-\nu(x))^2\right]\leq C_1\min_{h\in\mathcal H_n}\left\{\mathcal B_h(x)+\frac{\log(|\mathbb T_n|)}{|\mathbb T_n| |h|}\right\}+\frac{C_2}{|\mathbb T_n|},
\end{equation}
where $C_1,C_2>0$ do not depend on $n$ nor $x$ and 
$$\mathcal B_h(x)=\max_{h'\in\mathcal H_n}(K_{h'}\ast\nu(x)-K_h\ast K_{h'}\ast\nu(x))^2.$$
 \end{thm}

\begin{rk}
The form of the bias term $\mathcal B_h(x)$ in Inequality~\eqref{eq:ineg_oracle_ponctuelle} is very similar to the one obtained for pointwise adaptive kernel density estimation in \cite[Theorem 1]{Rebelles15b}. It can be replaced by an upper-bound, e.g. $\|K\|_1\|\nu-K_h\ast\nu\|_\infty$, coming from the Young's inequality, as in \cite[Theorem 2]{CL13}.
\end{rk}

Once the previous theorem is proved, an immediate corollary follows 
\begin{cor}\label{cor:adaptive_rates} Suppose that the assumptions of Theorem~\ref{thm:oracle} are verified and that $\nu\in\Sigma(\beta,L,D)$ where $\beta=(\beta_1,...,\beta_d)\in]0,+\infty[^d$ and $\Sigma(\beta,L,D)$ is the set of all $\beta$-H\"older densities on the open set $D\subseteq\RR^d$ i.e. the set of all functions $f:D\to\RR$ which admits, for all $j=1,...,d$, partial derivatives with respect to $x_j$ up to the order $\lfloor\beta_j\rfloor$ and verifies, for all $x=(x_1,...,x_d)\in D$, $x_j'\in\RR$ such that $(x_1,...,x_{j-1},x_j',x_{j+1},...,x_d)\in D$, 
$$\left|\frac{\partial^{\lfloor\beta_j\rfloor}f}{\partial x_j^{\lfloor\beta_j\rfloor}}(x_1,,...,x_{j-1},x_j',x_{j+1},...,x_d)-\frac{\partial^{\lfloor\beta_j\rfloor}f}{\partial x_j^{\lfloor\beta_j\rfloor}}(x)\right|\leq L|x_j-x_j'|^{\beta_j-\lfloor\beta_j\rfloor}.$$
Moreover, suppose that $K$ is a kernel of order $\ell\in\NN$ (with $\ell\geq\max_{j=1,...,d}\{\beta_j\}$) that is to say, $\int_\RR K(t)dt=1$, $\int_\RR x^jK(x)=0$ for all $j=1,...,\ell$ and that 
$$\int_\RR |x|^\ell |K(u)|<+\infty.$$
Suppose also that, for all $n$, there exists $h^*=(h_1^*,...,h_d^*)\in\mathcal H_n$ such that $h_j^*\asymp\left( \frac{\log(|\mathbb T_n|)}{|\mathbb T_n|}\right)^{\bar\beta/(\beta_j(2\bar\beta+d))}$ where $\bar\beta=d/(1/\beta_1+...+1/\beta_d)$ is the harmonic mean of $\beta$. 
Then, there exists a constant $C>0$ such that 
$$\sup_{x\in D}\EE\left[(\widehat\nu_{\widehat h}(x)-\nu(x))^2\right]\leq C\left(\frac{|\mathbb T_n|}{\log(|\mathbb T_n|)}\right)^{-2\bar\beta/(2\bar\beta+1)}.$$
\end{cor}
This corollary is a direct consequence of Theorem \ref{thm:oracle} and \cite{Tsybakov} Proposition 2.1. 

Remark that the assumption on the bandwidth collection is verified e.g. by $$\mathcal H_n=\left\{h_{\max} k^{-\alpha}, k=1,....,(|\mathbb T_n|h_{\max}/\log(|\mathbb T_n|))^{1/\alpha}\right\}^d$$ for all constant $h_{\max}>0$ and all $\alpha>1$.

\subsection{Estimation of the constant appearing in the variance}
\label{sec:slope_estim}
Due to the term $C(P,\mu)$ which is hardly estimable, the variance term $V(x,h)$ is not calculable in practice. We propose an algorithm to estimate it, based on slope estimation, as developed in \cite{AM09,LM16}. 

As suggested in \cite{LM16}, we take $b=2a$, hence in order to calculate the estimator, it is sufficient to have a \emph{good} value of $\kappa=aC(P,\mu)$. The following algorithm is inspired by the procedure described in \cite{AM09}, for model selection purposes. 
However, note that, in our case, both terms $A(x,h)$ and $V(x,h)$ depends on the constant $\kappa$, which is not the case in model selection contexts where only the penalty term depends on the constant. The selection of the grid of $\kappa$ is then different here. 

\begin{algo}
\begin{enumerate}
\item[]$\;$
\item[]$\;$
\item[1.] Initialization : set  
\begin{equation}\label{eq:kappam}
\kappa_m=\frac{|\mathbb T_n|}{\log(|\mathbb T_n|)}\max_{h,h'\in\mathcal H_n}{|h'|(\widehat\nu_{h'}(x)-K_{h}\ast\widehat\nu_{h'}(x))^2}\text{ and }\kappa_1=0.
\end{equation}
\item[2.] While $s\leq s_{\max}$
\begin{enumerate}
\item[(i)] Generate a sequence $(\kappa_j)_{1\leq j\leq m}$ such that 
 and $\kappa_j=\kappa_1+\frac{j-1}{m-1}(\kappa_m-\kappa_1)$, for all $j=1,...,m$. 
\item[(ii)] Calculate $\widehat h_j:=\widehat h(\kappa_j)$ as the minimizer of the criterion~\eqref{eq:critere} with $a=\kappa/ C(P,\mu)$ and $b=2a$.    
\item[(iii)] Set 
$$j_{jump}={\arg\max}_{j=1,...,m-1}{\left|\frac1{|\widehat h_j|}-\frac1{|\widehat h_{j+1}|}\right|},$$ 
$\kappa_1=\kappa_{j_{jump}}$,  $\kappa_1=\kappa_{j_{jump}+1}$ and $s=s+1$. 
\end{enumerate}
\item[3.] Return $\widehat h_{j_{jump}+1}$. 
\end{enumerate}
\end{algo}

The algorithm search the value of $\kappa$ for which the variance of the estimator increases significantly and select a slightly larger value. This allows to select an estimator with a reasonable variance. The same reasoning has given rise to the so-called \emph{dimension jump} method for model selection purposes (see \cite{AM09}). 
The chosen value for the initialization of $\kappa_m$ comes from the following reasoning. Setting $\kappa\geq\kappa_m$ as suggested: by definition, for all $h,h'\in\mathcal H_n$, 
$$\frac{|\mathbb T_n|h'}{\log(|\mathbb T_n|)}(\widehat\nu_{h'}(x)-K_{h}\ast\widehat\nu_{h'}(x))^2\leq \kappa$$
which implies that $A(x,h)=0$ for all $h\in\mathcal H_n$ and that the criterion~\eqref{eq:critere} will select the smaller bandwidth in $\mathcal H_n$. On the contrary, if $\kappa>\kappa_m$, let $h,h'\in\mathcal H_n$ for which the minimum is attained in~\eqref{eq:kappam}, we have 
$$\frac{|\mathbb T_n|}{\log(|\mathbb T_n|h')}(\widehat\nu_{h'}(x)-K_{h}\ast\widehat\nu_{h'}(x))^2> \kappa$$
which implies that $A(x,h)>0$ and the criterion may select a bandwidth which is not the smaller one. Hence, the values of $\kappa$ for which the criterion~\eqref{eq:critere} may select suitable values of the bandwidth can not be greater than $\kappa_m$. That is the reason why we consider an initial grid in the interval $[0,\kappa_m]$. 
However, the initial value $\kappa_m$ may be very large compared to the optimal value of $\kappa$. The loop 2. allows to search among small values of $\kappa$ while avoiding the choice of a too large $m$ which could be very expensive in terms of computation time. In practice $s_{\max}=2$ and $m=20$ seems to be a reasonable choice. 

\begin{rk}
The choice $b=2a$ is arbitrary since, in theory, any value $b$ such that $b\geq a$ might work as soon as $a$ sufficiently large. Other choices may be made such as $b=a$ which is the usual choice of Lepski's method, $b=a+\varepsilon$, with $\varepsilon>0$ or $b=(1+\varepsilon)a$. 
\end{rk}

\section{Simulations}
\label{sec:simus}

We shall now illustrate the previous results on simulated data coming from various models of bifurcating Markov chains.

\subsection{Bifurcating autoregressive processes}$\,$

Bifurcating autoregressive processes (BAR, for short)  were first
introduced by Cowan and Staudte \cite{CS86} in order to study the
data from cell division, where each individual in one generation
gives birth to two children in the next generation. This model has been widely studied over the last thirty years (see for e.g. \cite{BO16} and references therein). Recently, Bitseki and Olivier in \cite{BO16} have proposed an extension of BAR process initially introduced by Cowan and Staudte. Their model is defined as follows.

Let $X_{u} \in \RR$ be a quantitative data associated to
the cell $u \in \TT$, for example the growth rate of \textit{E.
Coli}. Then the quantities $X_{u0}$ and $X_{u1}$ associated to $u0$
and $u1$ the two children of $u$ are linked to $X_{u}$ through the
following autoregressive equations
\begin{equation}\label{eq BAR}
\begin{array}{ll}
\Ll(X_{\emptyset}) = \mu, \quad \text{and   for} \; u \in \TT, \quad
\left\{\begin{array}{ll} X_{u0} = f_0(X_{u}) + \varepsilon_{u0}, \\
\\ X_{u1} = f_1(X_{u})+ \varepsilon_{u1},
\end{array} \right.
\end{array}
\end{equation}
where $\mu$ is a distribution probability on $\RR$ and $f_0, f_1 :
\RR \leadsto \RR$. The noise $\big((\varepsilon_{u0},
\varepsilon_{u1}), u\in \TT\big)$ forms a sequence of independent
and identically distributed bivariate centered random variables with
common density $g_{\varepsilon}$ on $\RR \times \RR$. The process
$\left(X_{u}, u\in\TT\right)$ defined by \eqref{eq BAR} is a
bifurcating Markov chain with $\TT$-transition probability
\begin{equation*}
\Ttransition(x,dy,dz) =  g_{\varepsilon}(y-f_0(x),z-f_1(x)) dy dz.
\end{equation*}
Under some assumptions on $\mu$, $f_0$, $f_1$ and
$g_{\varepsilon}$, it has been shown in \cite{BO16} that the process $X$ satisfies all the good properties needed for our theoretical results (we refer to \cite{BO16} for more details).
We note that the previous model can be seen as an adaptation of nonlinear autoregressive model when the data have a binary tree structure. Furthermore, the original BAR process in \cite{CS86} is defined for linear link functions $f_0$ and $f_1$ with $f_0 = f_1$.

Now, for our numerical illustrations, we build a BAR process living in $S : = [0,1]$ as follows. First, we choose $X_{\emptyset}$ such that $\Ll(X_{\emptyset}) = Beta(2,2)$, where $Beta(2,2)$ is the standard Beta distribution with shape parameters $(2,2)$. Then for $u\in\TT$ and conditionally on $X_{u} = x$, we construct $X_{u0}$ and $X_{u1}$ independently in such a way that
$
\PP(X_{u0}\in dy,X_{u1}\in dz) = \Ttransition(x,y,z)dydz,
$
 where
$\Ttransition(x,\cdot,\cdot) : = \Pp(x,\cdot) \otimes
\Pp(x,\cdot)$ and
\begin{equation*} \label{BAR transition}
\Pp(x,y) : = (1 - x) \, \frac{y(1-y)^{2}}{B(2,3)} + x \,
\frac{y^{2}(1-y)}{B(3,2)} \, , \quad
x, y  \in [0,1]
\end{equation*}
with $B(\alpha, \beta)$ the normalizing constant of a standard Beta
distribution with shape parameters $\alpha$ and $\beta$. Now, one can prove that this process is stationary, it has an explicit invariant density, which is crucial to evaluate the quality of estimation of our method: this is a standard Beta distribution with shape parameters $(2,2)$. One can also prove that
$$
\EE\left[X_{u0}|X_{u}\right] = \EE\left[X_{u1}|X_{u}\right] = 1/5 X_{u} + 2/5,
$$ 
in such a way that the equations (\ref{eq BAR}) are satisfied with $f_{0}(x) = f_{1}(x) = 1/5 x + 2/5$ (for more details, we refer for e.g. to \cite{PCW02}). Now, it is no hard to verify that this process satisfies our required assumptions.

\begin{figure}
\includegraphics[width=0.45\textwidth]{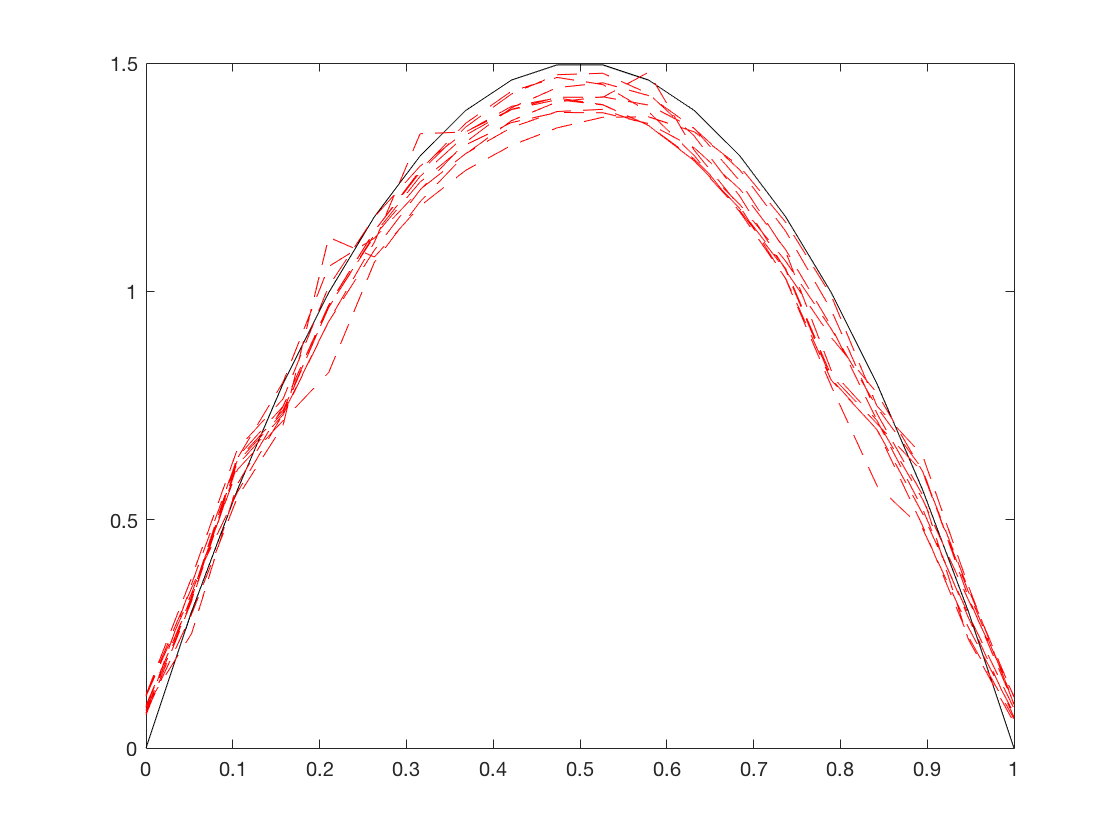}
\caption{\label{fig:estimsBAR} Plot of 10 estimators (red dashed lines) obtained from independent copies of $X$. The black solid line represents the function $\nu$ to estimate.  }
\end{figure}

We simulate the $n$ first generations of the process $X$, with $n=10$ (hence the size of $\mathbb T_n$ is $|\mathbb T_n|=2^{10}=1024$). We consider the Gaussian kernel $K(t)=\frac1{\sqrt{2\pi}}e^{-t^2/2}$.
The results are given in Figure~\ref{fig:estimsBAR}. 

\subsection{Estimation of splitting rate in a growth-fragmentation model}$\,$

We are now interested in growth-fragmentation models. Theses models describes (for e.g.) the evolution of cells which grow and divide randomly over time (see for e.g. \cite{DHKR1} and references therein). The model we are going to study is a simplification of the one studied in \cite{DHKR1}; it is defined as follows. Let $S$ be a subset of $[0,\infty)$ and let $B: S\mapsto [0,\infty)$ be a continuous function (the splitting rate). Each cell $u \in \TT$ grows exponentially with a common rate $\tau > 0$ and when it reaches a certain size $x$, it splits at rate $B(x)$, and gives birth to two offspring ($u0$ and $u1$) of size $x/2$. Next, this two offspring, $u0$ and $u1$, start a new life independently of each other. Clearly, the process $(X_{u},u\in\TT)$, where $X_{u}$ is the size of the cell $u$ at birth, is a BMC. It is proved in \cite{DHKR1} that the $\TT$-transition probability of this BMC is $\Ttransition(x,\cdot,\cdot) : = \Pp(x,\cdot) \otimes \Pp(x,\cdot)$, where the density of $\Pp(x,\cdot)$ is given by
\begin{equation*} \label{Pp_B}
\Pp (x, y) : = \frac{B(2y)}{\tau y} \exp \left( - \int_{x/2}^{y} \frac{B(2z)}{\tau z} dz \right) \mathds{1}_{\{y \geq x/2\}},
\end{equation*}
for $x \in S$ and $y \in S/2$. It can also be seen that the probability transition of the tagged-branch chain is $\Qq = \Pp$.

Doumic~\emph{et~al.} \cite{DHKR1}  have proved that $\Pp$ admits an invariant probability measure $\nu$ having a density, that we still denote by $\nu(\cdot)$, with respect to the Lebesgue measure. It is also known from \cite{DHKR1} that the rate function $B(\cdot)$  and the invariant density $\nu(\cdot)$ verify
$$
B(x) = \frac{\tau x}{2} \frac{\nu_B(x/2)}{\int_{x/2}^x \nu_B(z) dz}
$$ 
in such a way that a natural estimator for $B(x)$, based on the observation of $(X_{u}, u\in\TT_{n})$ is
\begin{equation*} \label{Bhat}
\widehat{B}_n(x) = \frac{\tau x}{2} \frac{\widehat{\nu}_h(x/2)}{\big( \tfrac{1}{|\TT_n|}\sum_{u \in \TT_n}{\bf 1}_{\{x/2 \leq X_u < x\}} \big)\vee \varpi_{n}},
\end{equation*} 
where $\widehat{\nu}_h$ is the kernel estimator defined in (\ref{eq:nuhat}) and $\varpi_{n}$ is a threshold which ideally go to $0$. Moreover, Bitseki~\emph{et~al.} \cite{BHO} have proved that under suitable assumptions on the splitting rate $B(\cdot)$, the process $X$ satisfies all the good properties needed for our theoretical results. For all the previous assertions, we refer to \cite{DHKR1,BHO} for more details. The strategy is then to use our results for the estimation of $\nu_{B}(x/2)$ for all $x\in S$. 

Now, for our numerical illustrations, we will work with the splitting rate using in \cite{BHO}. We choose $\tau = 2$, $S = (0,5)$ and for all $x\in S$, $B$ has the form 
\begin{equation*} \label{BHO_eq:Btrial}
B(x) = \frac{x}{x-5}  +  3 T\big(2 (x- \tfrac{7}{2})\big)
\end{equation*}
where 
$T(x)=(1+x){\bf 1}_{\{-1\leq x < 0\}}+(1-x){\bf 1}_{\{ 0 \leq x \leq 1 \}}$ is a tent shaped function. With this choice of $B$, the required assumptions for our theoretical results are satisfied (see \cite{BHO}).

\begin{figure}

\includegraphics[width=0.45\textwidth]{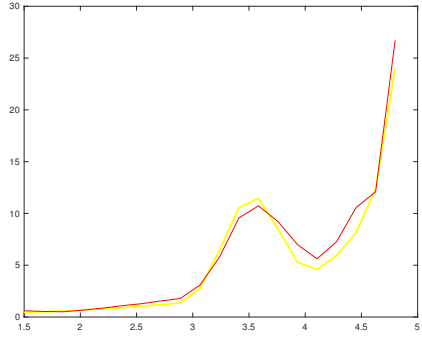}

\caption{\label{fig:estims-mcf} Plot of 1 estimator (red solid line) obtained from a copy of $X$. The yellow solid line represents the function $B$ to estimate.}
\end{figure}

We simulate the $n$ first generations of the process $X$, with $n=15$. The results are given in Figure~\ref{fig:estims-mcf}.

\section{Proofs}
\label{sec:proofs}
The proof rely on the lemma below, which is a Bernstein-type inequality. 
\begin{lem}\label{lem:bernstein} Let $(X_u, u\in\mathbb T_n)$ be a bifurcating Markov chain on $\RR^{d}$ with initial distribution $\mu$ and $\mathbb T$-transition probability $P$. Under the assumption of uniform geometric ergodicity, we have for all $\delta > 0$

\begin{multline}\label{eq:concentration2}
\PP\left(\frac{1}{|\TT_{n}|}\left|\sum_{u \in\TT_{n}}
\left(K_{h}*K_{h^{\prime}}\left(x - X_u \right) -
\EE_{\nu}\left[K_{h}*K_{h^{\prime}}\left(x - X_u
\right)\right]\right)\right|>\delta\right)
\\ \leq 2\exp\left(\frac{\delta c_{K,Q,\nu,M}c^{\prime}_{\rho}}{\frac{4c_{\rho,M}\|K\|_{1}\|K\|_{\infty}\delta}{3} + c_{K,Q,\nu,M}c^{\prime}_{\rho}}\right) \exp\left(-\frac{\delta^{2}|\TT_{n}||h^{\prime}|}{2\left(c_{K,Q,\nu,M} c^{\prime}_{\rho}
+ \frac{4c_{\rho,M}\|K\|_{1}\|K\|_{\infty}\delta}{3}\right)}\right),
\end{multline}

where 

\begin{multline*}
c_{\rho,M} = \frac{M(1+\rho)}{1-2\rho}, \, \, c^{\prime}_{\rho} = 3 + \frac{2}{1 - 2\rho}, \\ c_{K,Q,\nu,M} = 8 \max\{2\|K\|_{1}^{2}\|K\|_{2}^{2}(\|Q\|_{\infty,\infty}  + \|\nu\|_{\infty}); \max\{\|Q\|_{\infty,\infty} + \|\nu\|_{\infty}; M\|K\|_{1}\|K\|_{\infty}\}^{2}\}
\end{multline*}

We also have for all $\delta > 0$,

\begin{multline}\label{eq:concentration}
\PP\left(\frac{1}{|\TT_{n}|}\left|\sum_{u \in\TT_{n}}
\left(K_{h}\left(x - X_u\right) -
\EE_{\nu}\left[K_{h}\left(x - X_u
\right)\right]\right)\right| > \delta\right)
\\ \leq 2\exp\left(\frac{\delta c_{K,Q,\nu,M}^{\prime}c^{\prime}_{\rho}}{\frac{4c_{\rho,M} \|K\|_{\infty}\delta}{3} + c_{K,Q,\nu,M}^{\prime}c^{\prime}_{\rho}}\right) \exp\left(-\frac{\delta^{2}|\TT_{n}||h|}{2\left(c_{K,Q,\nu,M}^{\prime} c^{\prime}_{\rho}
+ \frac{4c_{\rho,M}\|K\|_{\infty}\delta}{3}\right)}\right)
\end{multline}
where $c_{K,Q,\nu,M}^{\prime} = 8 \max\{M \|K\|_{\infty}
; (\|\Qq\|_{\infty,\infty} + \|\nu\|_{\infty})\|K\|_{1};
(\|\Qq\|_{\infty,\infty} + \|\nu\|_{\infty})\|K\|_{2}^{2}\}$. 
\end{lem}

\begin{rk}
As mentioned above, these inequalities are more complete than those obtained in \cite{BHO}, since the deviation parameter $\delta$ does not depend on the size of the data. We stress that this fact is essential for our theoretical results.\end{rk}

\begin{proof}
We will do the proof of (\ref{eq:concentration2}). The proof of (\ref{eq:concentration}) follows the same lines.

Let $\lambda>0$ and $\delta>0$. By Chernoff inequality, we
have
\begin{multline*}
\PP\left(\frac{1}{|\TT_{n}|}\sum_{u \in\TT_{n}}
\left(K_{h}*K_{h^{\prime}}\left(x - X_u\right) -
\EE_{\nu}\left[K_{h}*K_{h^{\prime}}\left(x - X_u
\right)\right]\right) > \delta\right) \\ \leq
\exp\left(-\delta\lambda|\TT_{n}|\right)
\EE\left[\exp\left(\lambda\sum_{u\in\TT_{n}} g(X_{u})\right)\right],
\end{multline*}
where the function $g$ is defined by
\begin{equation*}
g(y) = K_{h}*K_{h^{\prime}}(x - y) -
\EE_{\nu}\left[K_{h}*K_{h^{\prime}}(x - X_{\emptyset}
)\right].
\end{equation*}

For all $u\in\GG_{n-1}$, we have, on the one hand
\begin{equation*}
|g(X_{u0}) + g(X_{u1}) - 2 \Qq g(X_{u})| \leq 2M(1 + \rho)\|K_{h}*K_{h^{\prime}}\|_{\infty}.
\end{equation*}

Using the Young's inequality, we have

\begin{equation*}
\|K_{h}*K_{h^{\prime}}\|_{\infty} \leq \|K\|_{1}\|K\|_{\infty}/|h'|
\end{equation*}

and therefore

\begin{equation*}
|g(X_{u0}) + g(X_{u1}) - 2\Qq g(X_{u}))| \leq 2 c_{\rho,M} \|K\|_{1}\|K\|_{\infty}/|h'|.
\end{equation*}

On the other hand, we have

\begin{equation*}
\EE\left[\left(g(X_{u0}) + g(X_{u1}) - 2\Qq g(X_{u})\right)^{2}|\Ff_{n-1}\right] \leq 4 \|Q\|_{\infty,\infty}\|K_{h}*K_{h^{\prime}}\|_{2}^{2}  \leq c_{K,Q,\nu,M} /|h^{\prime}|).
\end{equation*}

Now for all $\lambda \in
(0,3/(2c_{\rho,M}\|K\|_{1}\|K\|_{\infty}(|h^{\prime}|)^{-1}))$, we have from Bennett inequality

\begin{multline*}
\EE\left[\exp\left(\lambda\left(g\left(X_{u0}\right) +
g\left(X_{u1}\right) - 2\Qq g\left(X_{u}\right)\right)\right)
|\Ff_{n-1}\right]  \leq
\exp\left(\frac{c_{K,\Qq,\nu,M}(|h^{\prime}|)^{-1}\lambda^{2}}{2\left(1 - \frac{2c_{\rho,M}\|K\|_{1}\|K\|_{\infty}(|h^{\prime}|)^{-1}\lambda}{3}\right)}\right)
\end{multline*}

and using the Markov property,  this leads us to

\begin{multline*}
\EE\left[\exp\left(\lambda\sum_{u\in\TT_{n}} g(X_{u})\right)\right]
\leq \exp\left(\frac{c_{K,\Qq,\nu,M}(|h^{\prime}|)^{-1}\lambda^{2} |\GG_{n-1}|}{2\left(1 -
\frac{2c_{\rho,M}\|K\|_{1}\|K\|_{\infty}(|h^{\prime}|)^{-1}\lambda}{3}\right)}\right)
\\ \times \EE\left[\exp\left(\lambda\sum_{u\in\TT_{n-3}}
g(X_{u})\right) \times \exp\left(\lambda \sum_{u\in\GG_{n-2}}
I^{(2)}(X_{u})\right) \right. \\
\times \left. \prod_{u\in \GG_{n-2}}
\EE\left[\exp\left(J^{(2)}(X_{u},X_{u0},X_{u1})\right)\Big|\Ff_{n-2}\right]\right]
\end{multline*}

where 

\begin{align*}
I^{(2)}(X_{u}) &= (g + 2 \Qq g + 2^{2} \Qq^{2} g)(X_{u})& \\
J^{(2)}(X_{u},X_{u0},X_{u1}) &= (g + 2\Qq g)(X_{u0}) + (g + 2\Qq
g)(X_{u1}) - 2 (\Qq g + 2^{2}\Qq^{2}g)(X_{u}).&  
\end{align*}

Now for the second step, we will do the same thing with $J^{(2)}(X_{u},X_{u0},X_{u1})$ instead of $g(X_{u0}) + g(X_{u1}) - 2\Qq g(X_{u})$. For all $u\in \GG_{n-2}$ we have

\begin{equation*}
|J^{(2)}(X_{u},X_{u0},X_{u1})| \leq 2c_{\rho,M} \|K\|_{1} \|K\|_{\infty} (|h^{\prime}|)^{-1}.
\end{equation*}

We also have

\begin{equation*}
\EE[(J^{(2)}(X_{u},X_{u0},X_{u1}))^{2}|\Ff_{n-2}] \leq 4\Qq((g + 2\Qq g)^{2})(X_{u}).
\end{equation*}

Now we will control the right hand side of the previous inequality. On the one hand, we have for all $y\in \RR^{d}$,

\begin{equation*}
\Qq g(y) =  \int\int K_h\left(t\right) K_{h'}\left(x - z - t\right) \Qq(y,z) dtdz - \EE_{\nu}[K_{h}*K_{h^{\prime}}(x - X_{\emptyset}].
\end{equation*}

Using the change of variables

\begin{equation*}
t = hu \quad \text{and} \quad z = x - hu - h^{\prime}v,
\end{equation*}

where for $x =(x_1,...,x_d)^t,y=(y_1,...,y_d)^t\in\RR^{d}$, $xy$ denotes the vector $(x_1y_1,...,x_dy_d)^t$, we obtain

\begin{equation*}
\int\int K_h\left(t\right) K_{h'}\left(x-z-t\right) \Qq(y,z) dtdz| \leq \|Q\|_{\infty,\infty}.
\end{equation*}

In the same way we  prove that

\begin{equation*}
|\EE_{\nu}[K_{h}*K_{h^{\prime}}(x - X_{\emptyset}]| \leq \|\nu\|_{\infty}.
\end{equation*} 

This leads us to

\begin{equation*}
|\Qq g(y)| \leq \|Q\|_{\infty,\infty} + \|\nu\|_{\infty}
\end{equation*}

and therefore, for all $m\geq 1$

\begin{equation*}
|\Qq^{m} g(y)| \leq  \|Q\|_{\infty,\infty} + \|\nu\|_{\infty}.
\end{equation*}

On the other hand, uniform geometric ergodicity assumption implies that for all $m\geq 1$ and for all $y\in \RR^{d}$,

\begin{equation*}
\Qq^{m} g(y) \leq M \rho^{m} \|K_{h}*K_{h^{\prime}}\| \leq M \|K\|_{1} \|K\|_{\infty} \rho^{m} (|h^{\prime}|)^{-1},
\end{equation*} 

and therefore, for all $m\geq 1$ and for all $y\in \RR^{d}$, we have

\begin{equation*}
\Qq^{m} g(y) \leq \max\{\|Q\|_{\infty,\infty} + \|\nu\|_{\infty}, M \|K\|_{1} \|K\|_{\infty}\} \times \inf\{1,\rho^{m}(|h^{\prime}|)^{-1}\}.
\end{equation*}

For all $y\in \RR^{d}$, we also have

\begin{align*}
\Qq g^{2}(y) &\leq 2 \Qq(K_{h}*K_{h^{\prime}}) + 2 (\EE[K_{h}*K_{h^{\prime}}(x - X_{\emptyset})])^{2}& \\ &\leq 2 \|K\|_{1}^{2} \|K\|_{2}^{2} (\|Q\|_{\infty,\infty} + \|\nu\|_{\infty}),&
\end{align*}

in such a way that for all $u\in \GG_{n-2}$,

\begin{align*}
4\Qq((g + 2\Qq g)^{2})(X_{u}) \leq 8 \Qq g^{2}(X_{u}) + 8 \Qq((2 \Qq g)^{2})(X_{u}) \\ \leq c_{K,\Qq,\nu,M} ((|h^{\prime}|)^{-1} + 2 \inf\{1, \rho (|h^{\prime}|)^{-1}\}^{2})
\end{align*}

and thus, we obtain

\begin{equation*}
\EE[J^{(2)}(X_{u},X_{u0},X_{u1})|\Ff_{n-2}] \leq c_{K,\Qq,\nu,M} ((|h^{\prime}|)^{-1} + 2 \inf\{1, \rho (|h^{\prime}|)^{-1}\}^{2}).\end{equation*}

Once again, for all $\lambda \in
(0,3/(2c_{\rho,M}\|K\|_{1}\|K\|_{\infty}(|h^{\prime}|)^{-1}))$ and for all $u\in\GG_{n-2}$, we have from Bennett's inequality

\begin{equation*}
\EE[J^{(2)}(X_{u},X_{u0},X_{u1})|\Ff_{n-2}] \leq \exp\left(\frac{c_{K,\Qq,\nu,M}((|h^{\prime}|)^{-1} + 2 \inf\{1, \rho (|h^{\prime}|)^{-1}\}^{2})\lambda^{2}}{2\left(1 - \frac{2c_{\rho,M}\|K\|_{1}\|K\|_{\infty}(|h^{\prime}|)^{-1}\lambda}{3}\right)}\right).
\end{equation*}

It then follows that

\begin{multline*}
\EE\left[\exp\left(\lambda\sum_{u\in\TT_{n}} g(X_{u})\right)\right]
\leq \exp\left(\frac{c_{K,\Qq,\nu,M}(|h^{\prime}|)^{-1}\lambda^{2} |\GG_{n-1}|}{2\left(1 -
\frac{2c_{\rho,M}\|K\|_{1}\|K\|_{\infty}(|h^{\prime}|)^{-1}\lambda}{3}\right)}\right) \\ \times \exp\left(\frac{c_{K,\Qq,\nu,M}((|h^{\prime}|)^{-1} + 2 \inf\{1, \rho (|h^{\prime}|)^{-1}\}^{2})\lambda^{2}|\GG_{n-2}|}{2\left(1 - \frac{2c_{\rho,M}\|K\|_{1}\|K\|_{\infty}(|h^{\prime}|)^{-1}\lambda}{3}\right)}\right)\\ \times \EE\left[\exp\left(\lambda\sum_{u\in\TT_{n-4}}
g(X_{u})\right) \times \exp\left(\lambda \sum_{u\in\GG_{n-3}}
I^{(3)}(X_{u})\right) \right. \\
\times \left. \prod_{u\in \GG_{n-3}}
\EE\left[\exp\left(J^{(3)}(X_{u},X_{u0},X_{u1})\right)\Big|\Ff_{n-2}\right]\right]
\end{multline*}

where

\begin{align*}
I^{(3)}(X_{u}) &= (g + 2\Qq g + 2^{2}\Qq^{2} g + 2^{3} \Qq^{3}g)(X_{u});& \\ J^{(3)}(X_{u},X_{u0},X_{u1}) &= (g + 2\Qq g + 2^{2}\Qq^{2} g)(X_{u0}) + (g + 2\Qq g + 2^{2}\Qq^{2} g)(X_{u1}) - 2 (\Qq g + 2\Qq^{2} g + 2^{2} \Qq^{3}g)(X_{u}).&
\end{align*}

Now, iterating this method, we are led to

\begin{multline*}
\EE\left[\exp\left(\lambda\sum_{u\in\TT_{n}} g(X_{u})\right)\right]
\\ \leq \exp\left(\frac{\sum\limits_{m=1}^{n} c_{K,\Qq,\nu,M}\left((|h^{\prime}|)^{-1} + \left(\sum\limits_{l=1}^{m
-1} 2^{l} ((|h^{\prime}|)^{-1}\rho^{l} \wedge 1)\right)^{2}\right)|\GG_{n-m|}\lambda^{2}}{2\left(1
- \frac{2c_{\rho,M} \|K\|_{1} \|K\|_{\infty} (|h^{\prime}|)^{- 1}\lambda}{3}\right)}\right)
\\ \times \exp\left(c_{K,\Qq,\nu,M}\left((|h^{\prime}|)^{-1} + \sum_{m=1}^{n}
2^{m}\left(1 \wedge \rho^{m} (|h^{\prime}|)^{-1}\right)\right)\right).
\end{multline*}

Set $m^{*} = \lfloor\log |h^{\prime}|/\log\rho\rfloor$. Then we have

\begin{multline*}
\sum\limits_{m=1}^{n}\left((|h^{\prime}|)^{-1} + \left(\sum\limits_{l=1}^{m
-1}2^{l}((|h^{\prime}|)^{-1}\rho^{l}\wedge1)\right)^{2}\right)|\GG_{n-m|} \\
= \sum\limits_{m=m^{*} + 1}^{n}\left((|h^{\prime}|)^{-1} +
\left(\sum_{l=1}^{m^{*}} 2^{l} + \sum\limits_{l=m^{*}}^{m -1} 2^{l}
h^{-1} \rho^{l}\right)^{2}\right)|\GG_{n-m|}  +
\sum\limits_{m=1}^{m^{*}}\left((|h^{\prime}|)^{-1} +
\left(\sum\limits_{l=1}^{m -1}2^{l}\right)^{2}\right)|\GG_{n-m|}
\\ \leq \left(6 + (1 + \frac{1}{1 - 2\rho})^{2}\right) (|h^{\prime}|)^{-1}|\TT_{n}|.
\end{multline*}

We also have

\begin{equation*}
(|h^{\prime}|)^{-1} + \sum_{m=1}^{n} 2^{m}\left(1 \wedge
\rho^{m}(|h^{\prime}|)^{-1}\right) \leq c^{\prime}_{\rho} (|h^{\prime}|)^{-1}.
\end{equation*}

In view of the above, for all $\lambda \in
(0,3/(2c_{\rho,M} \|K\|_{1} \|K\|_{\infty}(|h^{\prime}|)^{-1}))$ we have

\begin{multline*}
\PP\left(\frac{1}{|\TT_{n}|}\sum_{u \in\TT_{n}}
\left(K_{h}*K_{h^{\prime}}(x - X_u) -
\EE_{\nu}\left[K_{h}*K_{h^{\prime}}(x - X_u)\right]\right)
> \delta\right) \\ \leq \exp\left(-\lambda \delta |\TT_{n}| +
\frac{c_{K,\Qq,\nu,M}c^{\prime}_{\rho} (|h^{\prime}|)^{-1}|\TT_{n}|\lambda^{2}}{2\left(1 -
\frac{2c_{\rho,M} \|K\|_{1} \|K\|_{\infty}(|h^{\prime}|)^{- 1}\lambda}{3}\right)}\right)
\times \exp\left(\lambda c_{K,\Qq,\nu,M}c^{\prime}_{\rho} (|h^{\prime}|)^{-1}\right).
\end{multline*}

Taking $\lambda = (\delta (|h^{\prime}|))/(c_{K,\Qq,\nu,M}c^{\prime}_{\rho} + (4c_{\rho,M} \|K\|_{1} \|K\|_{\infty}
\delta)/3)$, we obtain

\begin{multline*}
\PP\left(\frac{1}{|\TT_{n}|}\sum_{u \in\TT_{n}}
\left(K_{h}*K_{h^{\prime}}(x - X_u) -
\EE_{\nu}\left[K_{h}*K_{h^{\prime}}(x - X_u)\right]\right)
> \delta\right) \\ \leq \exp\left(\frac{\delta c_{K,\Qq,\nu,M}c^{\prime}_{\rho}}{\frac{4c_{\rho,M} \|K\|_{1} \|K\|_{\infty}\delta}{3} + c_{K,\Qq,\nu,M}c^{\prime}_{\rho}}\right) \exp\left(-\frac{\delta^{2}|\TT_{n}|(|h^{\prime}|)}{2\left(c_{K,\Qq,\nu,M} c^{\prime}_{\rho}
+ \frac{4c_{\rho,M} \|K\|_{1} \|K\|_{\infty}\delta}{3}\right)}\right).
\end{multline*}

The result follows since we can do the same thing for $-g$ instead
of $g$. Now, the proof of (\ref{eq:concentration}) follows the same lines and this ends the proof.
\end{proof}

\begin{proof}[Proof of Theorem~\ref{thm:oracle}]
We start from the following decomposition, true for all $h\in\mathcal H_n$
$$
(\widehat\nu_{\widehat h}(x)-\nu(x))^2\leq3(\widehat\nu_{\widehat h}(x)-K_h\ast\widehat\nu_{\widehat h}(x))^2+3(K_h\ast\widehat\nu_{\widehat h}(x)-\widehat\nu_h(x))^2+ 3(\widehat\nu_h(x)-\nu(x))^2. 
$$ 
Hence, since $K_h\ast\widehat\nu_{\widehat h}(x)=K_{\widehat h}\ast\widehat\nu_h(x)$, by definition of $A(x,h)$ and then by definition of $\widehat h$ and the fact that $a\leq b$
\begin{eqnarray}
(\widehat\nu_{\widehat h}(x)-\nu(x))^2&\leq&3((\widehat\nu_{\widehat h}(x)-K_h\ast\widehat\nu_{\widehat h}(x))^2-aV(x,\widehat h)+aV(x,\widehat h))\nonumber\\
&&+3((K_{\widehat h}\ast\widehat\nu_h(x)-\widehat\nu_h(x))^2-aV(x,h)+aV(x,h))+ 3(\widehat\nu_h(x)-\nu(x))^2\nonumber\\
&\leq& 3(A(x,h)+ bV(x,\widehat h))\nonumber\\
&&+3(A(x,\widehat h)+bV(x,h))+3(\widehat\nu_h(x)-\nu(x))^2\nonumber\\
&\leq & 6 (A(x,h)+bV(x,h)) +3(\widehat\nu_h(x)-\nu(x))^2\label{eq:oracle1}.
\end{eqnarray}

%
%

Now it remains to upper-bound $\mathbb E[A(x,h)]$.
We have
\begin{eqnarray*}
(\widehat\nu_{h'}(x)-K_h\ast\widehat\nu_{h'}(x))^2&\leq& 3(\widehat\nu_{h'}(x)-K_{h'}\ast\nu(x))^2+3(K_h\ast\widehat\nu_{h'}(x)-K_h\ast K_{h'}\ast\nu(x))^2\\
&&+3(K_{h'}\ast\nu(x)-K_h\ast K_{h'}\ast\nu(x))^2. 
\end{eqnarray*}

With a rough upper-bound of the $\max_{h\in\mathcal H_n}$ by the $\sum_{h\in\mathcal H_n}$ we get
\begin{eqnarray*}
\mathbb E[A(x,h)]&\leq &3\;\mathbb E\left[\max_{h'\in\mathcal H_n}\left((\widehat\nu_{h'}(x)-K_{h'}\ast\nu(x))^2-a\frac{V(x,h')}6\right)_+\right]\\
&&+\;3\;\mathbb E\left[\max_{h'\in\mathcal H_n}\left(K_h\ast\widehat\nu_{h'}(x)-K_h\ast K_{h'}\ast\nu(x))^2-a\frac{V(x,h')}6\right)_+\right]\\
&& +\; 3\;\max_{h'\in\mathcal H_n}(K_{h'}\ast\nu(x)-K_h\ast K_{h'}\ast\nu(x))^2\\
&\leq & 3 \sum_{h'\in\mathcal H_n}\mathbb E\left[\left(\left(\frac1{|\mathbb T_n|}\sum_{u\in\mathbb T_n}K_{h'}(x-X_u)-\mathbb E_\nu[K_{h'}(x-X_u)]\right)^2-a\frac{V(x,h')}6\right)_+\right]\\
&& \hspace{-2cm}+3 \sum_{h'\in\mathcal H_n}\mathbb E\left[\left(\left(\frac1{|\mathbb T_n|}\sum_{u\in\mathbb T_n}K_h\ast K_{h'}(x-X_u)-\mathbb E_\nu[K_h\ast K_{h'}(x-X_u)]\right)^2-a\frac{V(x,h')}6\right)_+\right]\\
&&+  3\;\max_{h'\in\mathcal H_n}(K_{h'}\ast\nu(x)-K_h\ast K_{h'}\ast\nu(x))^2\\
&\leq &\mathcal T_1+\mathcal T_2 +\mathcal B_h(x).
\end{eqnarray*}

We first give an upper-bound  for $\mathcal T_1$. Let $h'\in\mathcal H_n$ fixed, now remark that, by Lemma~\ref{lem:bernstein},
\begin{eqnarray*}
\mathcal T_1&=&\int_0^{+\infty}\mathbb P\left(\left(\left(\frac1{|\mathbb T_n|}\sum_{u\in\mathbb T_n}K_{h'}(x-X_u)-\mathbb E_\nu[K_{h'}(x-X_u)]\right)^2-a\frac{V(x,h')}6\right)_+\geq t\right)dt\\
&\leq&\int_0^{+\infty}\mathbb P\left(\left|\frac1{|\mathbb T_n|}\sum_{u\in\mathbb T_n}K_{h'}(x-X_u)-\mathbb E_\nu[K_{h'}(x-X_u)]\right|\geq \sqrt{t+a\frac{V(x,h')}6}\right)dt\\
&\leq&\int_{aV(x,h')/6}^{+\infty}\mathbb P\left(\left|\frac1{|\mathbb T_n|}\sum_{u\in\mathbb T_n}K_{h'}(x-X_u)-\mathbb E_\nu[K_{h'}(x-X_u)]\right|\geq \sqrt{u}\right)du\\
&\leq&\int_{aV(x,h')/6}^{+\infty}\exp\left(\frac{\sqrt uc_Kc_\rho'}{\frac{4c_\rho\|K\|_\infty\sqrt u}3+c_Kc_\rho'}\right)\exp\left(-\frac{u|\mathbb T_n||h|}{2(c_Kc'_\rho+(4/3)c_\rho\|K\|_\infty\sqrt u)}\right)du\\
&\leq & I_1+I_2,
\end{eqnarray*}
where 
\begin{eqnarray*}
I_1&=&\int_{aV(x,h')/6}^{aC(P,\mu)/6}\exp\left(\frac{\sqrt uc_Kc_\rho'}{\frac{4c_\rho\|K\|_\infty\sqrt u}3+c_Kc_\rho'}\right)\exp\left(-\frac{u|\mathbb T_n||h|}{2(c_Kc'_\rho+(4/3)c_\rho\|K\|_\infty\sqrt u)}\right)du\\
I_2&=&\int_{aC(P,\mu)/6}^{+\infty}\exp\left(\frac{\sqrt uc_Kc_\rho'}{\frac{4c_\rho\|K\|_\infty\sqrt u}3+c_Kc_\rho'}\right)\exp\left(-\frac{u|\mathbb T_n||h|}{2(c_Kc'_\rho+(4/3)c_\rho\|K\|_\infty\sqrt u)}\right)du,
\end{eqnarray*}
where we recall that, since for all $h'\in\mathcal H_n$, $|h'|\geq\frac{\log(|\mathbb T_n|}{|\mathbb T_n|}$, we have $V(x,h')\leq C(P,\mu)$. 

We first upper-bound $I_1$, 
\begin{eqnarray}
I_1&\leq& C_1\int_{aV(x,h')/6}^{aC(P,\mu)/6}\exp\left(-\frac{u|\mathbb T_n||h'|}{2(c_Kc'_\rho+(2\sqrt 2/3\sqrt 3)c_\rho\|K\|_\infty\sqrt{aC(P,\mu}))}\right)du\nonumber\\
     &\leq&\frac{C_1'}{|\mathbb T_n||h'|}\exp\left(-\sqrt{a}c_1^*|\mathbb T_n||h'|\right)\leq C_1'\log(|\mathbb T_n|)\exp\left(-\sqrt{a}c_1^*\log(|\mathbb T_n|)\right)=C_1'|\mathbb T_n|^{-\sqrt{a}c_1^*}\label{eq:controlT1}, 
\end{eqnarray}
with $C_1=\exp(\sqrt{aC(P,\mu)}/6)$, $C_1'=C_12(c_Kc'_\rho+(2\sqrt 2/3\sqrt 3)c_\rho\|K\|_\infty\sqrt{aC(P,\mu}))$, $c_1^*=C(P,\mu)/(12(c_Kc'_\rho+(2\sqrt 2/3\sqrt 3)c_\rho\|K\|_\infty\sqrt{C(P,\mu})))$, using the fact that $a\geq 1$ (remark that $c_1^*$ does not depend on $a$). 

We turn now to $I_2$, remark that the function $u\mapsto \sqrt uc_Kc_\rho'/(\frac{4c_\rho\|K\|_\infty\sqrt u}3+c_Kc_\rho')$ is non decreasing and converges to $3c_Kc_\rho'/(4c_\rho\|K\|_\infty)$ when $u\to\infty$, hence it is bounded by this quantity. We have then, using again $a\geq 1$, 
\begin{eqnarray}
I_2&\leq& C_2\int^{+\infty}_{aC(P,\mu)/6}\exp\left(-\frac{\sqrt{u}|\mathbb T_n||h'|}{2\left(\frac{\sqrt 6c_Kc'_\rho}{\sqrt{C(P,\mu)}}+(4/3)c_\rho\|K\|_\infty\right)}\right)du\nonumber\\
     &\leq&\left(\frac{C_2'}{|\mathbb T_n||h'|}+\frac{C_2''}{(|\mathbb T_n||h'|)^2}\right)\exp(-\sqrt ac_2^*|\mathbb T_n||h'|)\nonumber\\
     &\leq &C_2'\log(|\mathbb T_n|)|\mathbb T_n|^{-\sqrt ac_2^*}+C_2''\log^2(|\mathbb T_n|)|\mathbb T_n|^{-\sqrt ac_2^*}\label{eq:controlT2},
\end{eqnarray}
with 

\begin{align*}
C_2=\exp(3c_Kc_\rho'/(4c_\rho\|K\|_\infty), \quad C_2'=C_2\sqrt a(4c_Kc_\rho'+\frac{16}{3\sqrt 6}\sqrt{C(P,\mu)}c_\rho\|K\|_\infty),  \\
C_2''=8C_2\left(\sqrt 6c_Kc_\rho'/\sqrt{C(P,\mu)}+\frac43c_\rho\|K\|_\infty\right)^2 \quad \text{and} \quad c_2^*=\frac{\sqrt{C(P,\mu)}}{2\sqrt 6\left(\sqrt 6 c_Kc_\rho'/\sqrt{C(P,\mu)}+\frac43c_\rho\|K\|_\infty\right)}.
\end{align*}

Hence, gathering~\eqref{eq:controlT1} and~\eqref{eq:controlT2}, there exists $C'>0$ depending only on $C(P,\mu)$, $K$, $c_K$ and $\rho$, such that
\begin{eqnarray*}
\mathcal T_1&\leq& C\; \text{card}(\mathcal H_n)\log^2(|\mathbb T_n|)|\mathbb T_n|^{-\sqrt{a}c^*}\leq C' |\mathbb T_n|^{-1}
\end{eqnarray*}
with $C=3\max\{C_1',C_2',C_2''\}$, $c^*=\min\{c_1^*,c_2^*\}$ as soon as $a>4/(c^*)^2$.
\end{proof}

\section{Acknowledgement}
The authors want to thank Claire Lacour for her helpful advices on constant calibration.

\end{document}